\documentstyle[12pt, page1]{article}
\newcommand{\sect}[1]{\section{#1}\setcounter{equation}{0}}

\font\mbn=msbm10 scaled \magstep1
\font\mbs=msbm7 scaled \magstep1
\font\mbss=msbm5 scaled \magstep1
\newfam\mbff
\textfont\mbff=\mbn
\scriptfont\mbff=\mbs
\scriptscriptfont\mbff=\mbss
\def\mbf{\fam\mbff}
\def\Re{{\mbf R}}

\def\Co{{\mbf C}}

\newtheorem{Th}{Theorem}[section]

\newtheorem{C}[Th]{Corollary}

\newtheorem{R}[Th]{Remark}

\author{Alexander Brudnyi\thanks
{Research supported in part by NSERC.
  }
\\ Department of 
Mathematics and Statistics\\ University of Calgary, Calgary\\ Canada} 
\title{ THE DISTRIBUTION OF VALUES OF ANALYTIC FUNCTIONS ON 
CONVEX SETS
}
\date{September 15, 2000}
\begin{document}
\maketitle
\begin{abstract}
{Proceeding the study of local properties of analytic functions started in 
[Br] we prove new dimensionless inequalities for such functions in terms
of their Chebyshev degree. As a consequence we obtain the reverse H\"{o}lder
inequalities for analytic functions with absolute (i.e., independent of
dimension) constants. For polynomials such inequalities were recently proved
by Bobkov who sharpened and generalized the previous Bourgain result and
by Sodin and Volberg.}
\end{abstract}
\sect{\hspace*{-1em}. Introduction.}
{\bf 1.1.} 
In order to formulate the main result we first recall the definition of
the {\em Chebyshev degree} for analytic functions.

Let $B_{c}(0,1)\subset B_{c}(0,r)\subset\Co^{n}$ be the pair of open
complex Euclidean balls of radii 1 and $r$ centered at 0. Denote by
${\cal O}_{r}$ the set of holomorphic functions defined on $B_{c}(0,r)$.
Our definition is motivated by the following result (see [Br, Th. 1.1]).
\begin{Th}\label{te0}
Let $f\in {\cal O}_{r}$, $r>1$, and $I$ be a real interval situated in
$B_{c}(0,1)$. (Hereafter we identify $\Co^{n}$ with $\Re^{2n}$.) There
is a constant $d=d(f,r)>0$ such that for any $I$ and any measurable subset
$\omega\subset I$
\begin{equation}\label{eq1}
\sup_{I}|f|\leq\left(\frac{4|I|}{|\omega|}\right)^{d}\sup_{\omega}|f|\ .
\end{equation}
\end{Th}
The optimal constant in this inequality is called the {\em Chebyshev degree}
of $f\in {\cal O}_{r}$ in $B_{c}(0,1)$ and is denoted by $d_{f}(r)$. According
to the classical Remez inequality $d_{f}(r)$ does not exceed the (total)
degree of $f$, provided $f$ is a polynomial. Example 1.14 in [Br] shows
that even in this case it can be essentially smaller than degree.

Further, let us recall that a function $h:\Re^{k}\longrightarrow\Re_{+}$ is
${\em log-concave}$ if its support
$$
K=\{x\in\Re^{k}:\ h(x)>0\}
$$
is convex and $\log h$ is a concave function on $K$.\\
Examples of log-concave functions. (a) 
A nonnegative  function which is concave on a convex body in $\Re^{k}$ and is
zero outside. In particular, the indicator function of a convex body is 
log-concave.\\
(b) Many density functions of statistics, e.g., $e^{-|x|^{\alpha}}$, 
$\alpha>0$.\\
(c) $\max(0,l^{k})$ where $l$ is linear and $k$ is a positive integer.

Let now $\mu_{h}$ be a measure on $\Re^{k}$ with density $h$. For a convex
body $V\subset\Re^{k}$ set
$$
|V|:=\mu_{h}(V),\ \ \  f_{V}:=
\exp\left(\frac{1}{|V|}\int_{V}\log|f|d\mu_{h}\right)\ .
$$
Clearly,
$$
f_{V}\leq\frac{1}{|V|}\int_{V}|f|d\mu_{h}\ .
$$
In the formulation of the main result we assume without loss of generality
that
$$
\mu_{h}(V)=1\ .
$$
\begin{Th}\label{te1}
Let $\Re^{k}\subset\Co^{n}(\cong\Re^{2n})$ be a $k$-dimensional affine 
subspace,
$V\subset B_{c}(0,1)\cap\Re^{k}$ be a $k$-dimensional convex body and
$h:\Re^{k}\longrightarrow\Re_{+}$ be a log-concave function supported on $V$.
There are absolute constants 
$c,C>0$ (i.e. independent of dimensions $k, n$) such that for every $r>1$ and
$f\in {\cal O}_{r}$
\begin{equation}\label{eq2}
\begin{array}{lr}
\displaystyle
(1)\ \ \ \mu_{h}\{x\in V:\ |f(x)|>\alpha f_{V}\}\leq
C\exp(-c\alpha^{1/d_{f}(r)})\\
\\
{\rm and}\\
\displaystyle
\\
(2)\ \ \ \mu_{h}\{x\in V:\  |f(x)|\leq\alpha f_{V}\}\leq 
C(c\alpha)^{1/d_{f}(r)}\left(\log\alpha\right)^{1/2}, \alpha\leq e^{-1} \ .
\end{array}
\end{equation}
\end{Th}
\begin{C}\label{c1}
Under the assumptions of Theorem \ref{te1}
$$
\frac{1}{|V|}\int_{V}|f|^{p}d\mu_{h}\leq (cpd_{f}(r))^{pd_{f}(r)}(f_{V})^{p}
\leq (cpd_{f}(r))^{pd_{f}(r)}
\left(\frac{1}{|V|}\int_{V}|f|d\mu_{h}\right)^{p}\ \ \ 
(p>1)
$$ 
with an absoulte constant $c>0$.
\end{C}

In particular, if $f_{V}\leq 1$, then the Orlicz norm of $f$ defined
by the Orlicz function $\exp(1/d_{f}(r))-1$ on $(V,d\mu_{h})$ is bounded
by an absolute constant.
\begin{C}\label{c2}
Under assumptions of Theorem \ref{te1}
$$
\frac{1}{|V|}\int_{V}\left|\log |f|-C_{V}(f)\right|d\mu_{h}\leq
Cd_{f}(r)\ .
$$
Here $C>0$ is an absoulte constant and $C_{V}(f):=
\frac{1}{|V|}\int_{V}\log |f|d\mu_{h}$.
\end{C}
A similar result for analytic functions compairing $\log |f|$
with $\sup_{V}\log |f|$ was obtained in [Br]. In this case the constant
in the inequality is equivalent to $\log k$ (for $k=dim V$).
In the case of polynomials Corollary \ref{c1} implies the fundamental
Bourgain inequality [B] (with $h\equiv 1$) and its generalizations proved
by Bobkov [Bo] and by the author in [Br, Th.1.11]. 
Inequality $(2)$ implies a similar inequality for 
plurisubharmonic functions on $\Co^{n}$ of a logarithmic growth
recently proved in [SV]. This paper contains also 
Corollary \ref{c2} for such functions.

{\bf 1.2.} 
As in the above cited papers [Bo] and [SV] our main tool is a remarkable
result of Kannan, Lov\'{a}sz and Simonovits ([KLS, Cor. 2.2]) which
reduces estimation of a multidimensional integral to the corresponding
one-dimensional ones. Using this we establish the following basic
inequality which gives Theorem \ref{te1} as a simple consequence.
\begin{Th}\label{te2}
Let $f$, $V$, $r>1$ and $\mu_{h}$ be as in Theorem \ref{te1}. Then
$$
\left(\frac{1}{|V|}\int_{V}|f|^{m}d\mu_{h}\right)^{n}
\left(\frac{1}{|V|}\int_{V}|f|^{-p}d\mu_{h}\right)^{q}\leq
(2e)^{n+q}4^{(mn+pq)d_{f}(r)}\frac{\Gamma(md_{f}(r)+1)^{n}}{(1-pd_{f}(r))^{q}}
$$
provided $m,n,p,q>0$ satisfy
$$
mn=pq,\ \ \ p<\frac{1}{d_{f}(r)}\ .
$$
Here, as usual, $\Gamma(x):=\int_{0}^{\infty}t^{x-1}e^{-t}dt$.
\end{Th}
{\em Aknowledgement.} I would like to thank Prof. M.Sodin for sending me 
the preprint of his and Prof. A.Volberg unpublished paper [SV].
\sect{\hspace*{-1em}. Proof of Theorem \ref{te2}.}
{\bf 2.1.} In this section we collect the results used in the proof of the 
theorem. First introduce the following definition (see [KLS]).

By an {\em exponential needle} we mean a segment $I=[a,b]$ in $\Re^{n}$, 
together with a real constant $\gamma$. If $(E,\gamma)$ is an exponential 
needle and $f$ is an integrable function defined on $I$, then we set
$$
\int_{E}f=\int_{0}^{|b-a|}f(a+tu)e^{\gamma t}dt,
$$
where $u=(1/|b-a|)(b-a)$.
\begin{Th}\label{KLS}{\bf [KLS]}
Let $f_{1}, f_{2},f_{3},f_{4}$ be four nonnegative continuous
functions defined on $\Re^{n}$, and $\alpha,\beta>0$. Then the following
are equivalent:

(a) For every log-concave function $F$ defined on $\Re^{n}$ with compact
support,
$$
\left(\int_{\Re^{n}}F(t)f_{1}(t)dt\right)^{\alpha}
\left(\int_{\Re^{n}}F(t)f_{2}(t)dt\right)^{\beta}\leq\\
\left(\int_{\Re^{n}}F(t)f_{3}(t)dt\right)^{\alpha}
\left(\int_{\Re^{n}}F(t)f_{4}(t)dt\right)^{\beta}\ .
$$

(b) For every exponential needle $E$ 
$$
\left(\int_{E}f_{1}\right)^{\alpha}\left(\int_{E}f_{2}\right)^{\beta}
\leq \left(\int_{E}f_{3}\right)^{\alpha}\left(\int_{E}f_{4}\right)^{\beta}\ .
$$
\end{Th}
\begin{R}\label{re1} 
{\rm The above theorem is also valid for nonnegative
$f_{1},f_{2},f_{3},f_{4}$ such that $f_{1},f_{2}$ are the limits of 
monotone increasing sequences of continuous functions defined on $\Re^{n}$ and
$f_{3},f_{4}$ are the limits of monotone decreasing sequences of 
continuous functions defined on $\Re^{n}$ (see Remark 2.3 in [KLS]). 
In particular, we can apply this theorem in the case if $K$ is a closed 
convex body, $f_{1},f_{2}$ are nonnegative continuous functions defined on 
$K$ which are 0 outside $K$ and $f_{3},f_{4}$ are nonnegative functions
which are constant on $K$ and 0 outside.}
\end{R}

We also use the following distributional inequality that follows directly
from inequality (\ref{eq1}) (see [Br]). \\
Let $I\subset B_{c}(0,1)
\subset\Co^{n}$ be a real segment and $f\in {\cal O}_{r}$.
For the distribution function $D_{f_{I}}(t):=|\{x\in I\ :\ |f(x)|\leq t\}|$
(with respect to the usual Lebesgue measure on $I$)
define $(f_{I})_{*}(t)=\inf\{s\ :\ D_{f_{I}}(s)\geq t\}$. Then
\begin{equation}\label{dist}
(f_{I})_{*}(t)\geq\left(\frac{t}{4|I|}\right)^{d_{f}(r)}
\sup_{V}|f|\ .
\end{equation}
{\bf 2.2.} {\bf Proof of Theorem \ref{te2}.}
Let $f\in {\cal O}_{r}$ and $I\subset B_{c}(0,1)$ be a real interval.
Then the functions $f_{\epsilon}:=(|f|+\epsilon)|_{I}$, $\epsilon>0$, and
$f_{\epsilon,a,b}(t)=f_{\epsilon}(at+b)$, $t\in I, a,b\in\Re$, also 
satisfy inequality
(\ref{eq1}). We must apply the KLS theorem to functions
$f_{1}:=(|f|+\epsilon)^{m}$, $f_{2}:=(|f|+\epsilon)^{-p}$
(continuous on $V$) and $f_{3}:=2e\cdot 4^{md_{f}(r)}\Gamma(md_{f}(r)+1)$,
$f_{4}:=2e\cdot 4^{pd_{f}(r)}/(1-pd_{f}(r))$ on $V$ and 0 outside $V$
and then take the limit when $\epsilon\to 0$. To avoid abuse of notation
and because our estimates below do not depend on $\epsilon$ we may
assume without loss of generality that $|f|$ itself has no zeros
on $B_{c}(0,1)$.

According to the KLS theorem and Remark \ref{re1} the theorem follows from
the inequality
$$
\left(\int_{E}|f|^{m}\right)^{n}
\left(\int_{E}|f|^{-p}\right)^{q}\leq
(2e)^{n+q}4^{(mn+pq)d_{f}(r)}\frac{\Gamma(md_{f}(r)+1)^{n}}{(1-pd_{f}(r))^{q}}
\left(\int_{E}1
\right)^{n+q}
$$
for an exponential needle $E\subset V$.
Making an affine change of variables in the above integrals we
reduce the problem to the following inequality
$$
\begin{array}{l}
\displaystyle
\left(\int_{0}^{s}|\tilde f (x)|^{m}e^{-x}dx\right)^{n}
\left(\int_{0}^{s}|\tilde f (x)|^{-p}e^{-x}dx\right)^{q}\leq\\
\\
\displaystyle
(2e)^{n+q}4^{(mn+pq)d_{f}(r)}\frac{\Gamma(md_{f}(r)+1)^{n}}{(1-pd_{f}(r))^{q}}
(1-e^{-s})^{n+q}\ .
\end{array}
$$
Here $\tilde f$ is a function obtained from $f$ by this change of variables. 
As we already mentioned  $\tilde f$ satisfies (\ref{eq1}). Below we denote
$||\tilde f||_{I}:=\sup_{I}|\tilde f|$.

First, let $0\leq s\leq 1$. Then
$$
\begin{array}{c}
\displaystyle
\left(\int_{0}^{s}|\tilde f(x)|^{m}e^{-x}dx\right)^{n}
\left(\int_{0}^{s}|\tilde f(x)|^{-p}e^{-x}dx\right)^{q}\leq\!
\left(\int_{0}^{s}\left(\frac{|\tilde f(x)|}{||\tilde f||_{[0,s]}}\right)^{m}
\!dx\right)^{n}
\left(\int_{0}^{s}\left(\frac{||\tilde f||_{[0,s]}}{|\tilde f(x)|}
\right)^{p}\!dx
\right)^{q}\\
\displaystyle
\leq s^{n}
\left(\int_{0}^{s}\left(\frac{||\tilde f||_{[0,s]}}
{\tilde f_{*}(t)}\right)^{-p}dt\right)^{q}
\leq s^{n}
\left(s\int_{0}^{1}\left(\frac{4}{t}\right)^{pd_{f}(r)}dt\right)^{q}
\leq 4^{pqd_{f}(r)}s^{n+q}\left(\frac{1}{1-pd_{f}(r)}\right)^{q}\leq\\
\displaystyle
4^{pqd_{f}(r)}(2(1-e^{-s}))^{n+q}\left(\frac{1}{1-pd_{f}(r)}\right)^{q}\ .
\end{array}
$$
Here we applied inequality (\ref{dist}) to the lower distribution function 
$\tilde f_{*}$
of $\tilde f$ and used the inequality $1-e^{-s}>s/2$ for $0<s\leq 1$. Observe
that the obtaining constant is even less than the required one.

Assume now that $s>1$. We estimate each of two factors of the given
expression. Without loss of generality we may assume that $s$ is an integer.
Then
$$
\begin{array}{c}
\displaystyle
\int_{0}^{s}|\tilde f(x)|^{m}e^{-x}dx=\sum_{i=0}^{s-1}
\int_{i}^{i+1}|\tilde f(x)|^{m}e^{-x}dx\leq
\sum_{i=0}^{s-1}\left(\int_{i}^{i+1}
|\tilde f(x)|^{m}dx\right)e^{-i}\leq\\
\displaystyle
\sum_{i=0}^{s-1}\left(\int_{i}^{i+1}\left(
\frac{|\tilde f(x)|}{||\tilde f||_{[i,i+1]}}
||\tilde f||_{[i,i+1]}\right)^{m}dx\right)e^{-i}\leq
\sum_{i=0}^{s-1}||\tilde f||_{[0,i+1]}^{m}e^{-i}\leq\\
\displaystyle
\sum_{i=0}^{\infty}(4(i+1))^{md_{f}(r)}e^{-i}||\tilde f||_{[0,1]}^{m}\leq
4^{md_{f}(r)}e\int_{0}^{\infty}x^{md_{f}(r)}e^{-x}dx 
||\tilde f||_{[0,1]}^{m}=\\
\displaystyle
4^{md_{f}(r)}e\Gamma(md_{f}(r)+1)||\tilde f||_{[0,1]}^{m}\ .
\end{array}
$$
We used here inequality (\ref{eq1}) to estimate $\sup_{[0,i+1]}|\tilde f|$ by
$\sup_{[0,1]}|\tilde f|$.
Similarly,
$$
\begin{array}{c}
\displaystyle
\int_{0}^{s}|\tilde f(x)|^{-p}e^{-x}dx\leq\sum_{i=0}^{s-1}\left(\int_{i}^{i+1}
|\tilde f(x)|^{-p}dx\right)e^{-i}\leq\\
\displaystyle
\sum_{i=0}^{s-1}\left(\int_{0}^{i+1}\left(\frac{||\tilde f||_{[0,i+1]}}
{|\tilde f(x)|}\right)^{p}\frac{1}{||\tilde f||_{[0,i+1]}^{p}}dx\right)e^{-i}
\leq\\
\displaystyle
\sum_{i=0}^{s-1}\left(\int_{0}^{i+1}\left(\frac{||\tilde f||_{[0,i+1]}}
{\tilde f_{*}(t)}\right)^{p}\frac{1}{||\tilde f||_{[0,1]}^{p}}dt\right)e^{-i}
\leq\\
\displaystyle
\left(\sum_{0}^{s-1}\int_{0}^{i+1}\left(\frac{4(i+1)}{t}\right)
^{pd_{f}(r)}\frac{1}
{||\tilde f||_{[0,1]}^{p}}dt\right)e^{-i}\leq\\
\displaystyle
\sum_{i=0}^{s-1}\frac{4^{pd_{f}(r)}(i+1)}{1-pd_{f}(r)}\frac{1}{
||\tilde f||_{[0,1]}^{p}}e^{-i}\leq\frac{e 4^{pd_{f}(r)}}{1-pd_{f}(r)}
\frac{1}{||\tilde f||_{[0,1]}^{p}}\ .
\end{array}
$$
Using that $pq=mn$ and
$1-e^{-s}\geq 1/2$ for $s\geq 1$ we get from these inequalities
$$
\begin{array}{c}
\displaystyle
\left(\int_{0}^{s}|\tilde f(x)|^{m}e^{-x}dx\right)^{n}
\left(\int_{0}^{s}|\tilde f(x)|^{-p}e^{-x}dx\right)^{q}\leq\\
\displaystyle
4^{(mn+pq)d_{f}(r)}\cdot (2e)^{n+q}
\frac{\Gamma(md_{f}(r)+1)^{n}}{(1-pd_{f}(r))^{q}}(1-e^{-s})^{n+q}\ .
\end{array}
$$

This completes the proof of the theorem. \ \ \ \ $\Box$
\sect{\hspace*{-1em}. Proof of Theorem \ref{te1} and Corollaries.}
{\bf Proof of Theorem \ref{te1}.}
(1) We apply Theorem \ref{te2} to $g:=|f|^{1/d_{f}(r)}$ with
$n=1$, $p=1/2$, $q=2m$ and $m$ a positive integer.
Assume without loss of generality that $g_{V}=1$ and set
$E_{\alpha}:=\{x\in V : g(x)>\alpha\}$, $|E_{\alpha}|:=\mu_{h}(E_{\alpha})$.
Then from Theorem \ref{te2} we obtain
$$
\alpha^{m}|E_{\alpha}|\left(\int_{V}g^{-1/2}d\mu_{h}\right)^{2m}\leq
\left(\int_{V}g^{m}d\mu_{h}\right)\left(\int_{V}g^{-1/2}
d\mu_{h}\right)^{2m}\leq
4^{2m}(2e)^{2m+1}2^{2m}(m!)
$$
which is equivalent to
\begin{equation}\label{divi}
\begin{array}{c}
\displaystyle
\alpha^{m}|E_{\alpha}|\leq\frac{2^{8m+1}e^{2m+1}(m!)}{(\int_{V}g^{-1/2}
d\mu_{h})^{2m}}\leq
2^{8m+1}e^{2m+1}(m!)\exp\left(-2m\log\left(\int_{V}g^{-1/2}d\mu_{h}\right
)\right)\\
\\
\displaystyle
\leq 2^{8m+1}e^{2m+1}(m!)(g_{V})^{m}=2^{8m+1}e^{2m+1}(m!)\leq e^{10m}(m!)\ .
\end{array}
\end{equation}
We used here the Jensen inequality
$$
\int_{V}g^{-1/2}d\mu_{h}\geq\exp\left(\frac{-1}{2}
\int_{V}\log g\ d\mu_{h}\right)\ .
$$
Since $|V|=1$, we also have
$$
|E_{\alpha}|\leq 1\ .
$$
Dividing both sides of (\ref{divi}) by  $e^{11m}(m!)$
and summing by $m$ from $0$ to $\infty$ we get
$$
\exp(\alpha/e^{11})|E_{\alpha}|\leq 2,
$$
or
$$
|E_{\alpha}|\leq 2\exp(-\alpha/e^{11}) .
$$
Since $g:=|f|^{1/d_{f}(r)}$, the required inequality follows from here.

This proves part (1).

(2) Recall that
$C_{V}(f):=\frac{1}{|V|}\int_{V}\log|f| d\mu_{h}$ .
We will estimate the measure $|F_{\gamma}|:=\mu_{h}(F_{\gamma})$ of the 
set $F_{\gamma}:=\{x\in V\  : |\log|f|- C_{V}(f)|\geq\gamma\}$,
$\gamma\geq 1$.
We apply Theorem \ref{te2} to $f$ with $m=p=(1-1/\gamma)/d_{f}(r)$, $n=q=1$.
Then we have
$$
\begin{array}{c}
e^{\frac{2\gamma (1-1/\gamma)}{d_{f}(r)}}|F_{\gamma}|^{2}\leq
\left(\int_{V}e^{\frac{1-1/\gamma}{d_{f}(r)}(\log|f|-C_{V}(f))}d\mu_{h}
\right)
\left(\int_{V}e^{\frac{-(1-1/\gamma)}{d_{f}(r)}
(\log|f|-C_{V}(f))}d\mu_{h}\right)\\
\\
\displaystyle
\leq 4^{2(1-1/\gamma)}(2e)^{2}\frac{\Gamma(2-1/\gamma)}{1-(1-1/\gamma)}\leq
2^{6}e^{2}\gamma
\end{array}
$$
Hence
\begin{equation}\label{dibmo}
|F_{\gamma}|\leq 2^{3}e^{1+1/d_{f}(r)}e^{-\gamma/d_{f}(r)}\gamma^{1/2}\ .
\end{equation}
This, in particular, gives an estimate of $\mu_{h}\{x\in V\ : 
\log |f|-C_{V}(f)\leq -\gamma\}$ which, in turn, gives the 
required result
$$
\mu_{h}\{x\in V\ : |f(x)|\leq\alpha f_{V}\}\leq 2^{3}e
(e\alpha)^{1/d_{f}(r)}\left(\log\alpha\right)^{1/2}, \alpha\leq e^{-1}\ .
$$
with $\alpha=e^{-\gamma}$.

The proof of Theorem \ref{te2} is complete.\ \ \ $\Box$

{\bf Proofs of Corollaries.}
Corollary \ref{c1} follows directly by integration of inequality
(1) of Theorem \ref{te1} and 
Corollary \ref{c2} is a simple consequence of inequality (\ref{dibmo}).
\ \ \ $\Box$
\sect{\hspace*{-1em}. Concluding Remarks.}
The basic point of the above proof is, together with the KLS theorem, the
inequality (\ref{eq1}). This type of inequalities hold for a more general
class of functions introduced as follows.

Let $B_{c}(0,1)\subset B_{c}(0,r_{1})\subset B_{c}(0,r)\subset\Co^{n}$, 
$1<r_{1}<r$, be open complex Euclidean balls. Further, let
$l\subset\Co^{n}$ be a complex straight line which intersects all
the balls and $l_{1}:=l\cap B_{c}(0,1),\ l_{r_{1}}:=l\cap 
B_{c}(0,r_{1})$. Let $f$ be a plurisubharmonic function defined on
$B_{c}(0,r)$. We set
$$
b_{f}(l,r_{1})=\sup_{l_{r_{1}}}f-\sup_{l_{1}}f\ .
$$
Then the Bernstein index $b_{f}(r_{1})$ is defined as supremum of
$b_{f}(l,r_{1})$ taken over all lines $l$ (see also Definition 2.1 in
[Br]). The argument for the proof of Th.1.1 in [Br] leads to the
following result.
\begin{Th}\label{conclud}
Let $f$ be as above and $b_{f}(r_{1})<\infty$. Then 
the function $F:=e^{f}$ satisfies inequality (\ref{eq1}) with $cb_{f}(r_{1})$
instead of $d_{f}(r)$ where $c\geq 1$ depends on $r_{1}$ and $r$ only.
\end{Th}
\begin{C}\label{colast}
All results of the presents paper are valid for $e^{f}$ with
$cb_{f}(r)$ instead of $d_{f}(r)$.
\end{C}
{\bf Examples.} (a) Let $f\in {\cal O}_{r}$ then we proved in [Br]
that $b_{\log |f|}((1+r)/2)<\infty$. This motivates our
definition of the Chebyshev degree. 

Note also that $d_{f}(r)$ can be 
estimated by the general valency of $f$ defined as maximum of valency
of $f$ restricted to each complex disk $l_{(1+r)/2}$.\\
(b) Let $f_{1},...,f_{k}\in(\Co^{n})^{*}$ be complex linear
functionals. A quasipolynomial with the spectrum $f_{1},...,f_{k}$
is a finite sum $q(z)=\sum_{i=1}^{k}p_{i}(z)e^{f_{i}(z)}$ where
$p_{i}\in\Co [z_{1},...,z_{n}]$ are holomorphic polynomials.
Expression $m=\sum_{i=1}^{k}(1+deg\ p_{i})$ is said to be degree of $q$.
We set
$$a(q):=\max_{z\in B_{c}(0,1), 1\leq i\leq k}|f_{i}(z)|. $$
From the results of [Br1] (see also [Br, Prop. 1.4]) it follows.

{\em All results of the present
paper are valid for $q|_{V}$ ($V\subset\Re^{n}$ is a convex body)
with $c_{1}+mc_{2}+c_{3}a(q)diam(V)$ instead of
$d_{f}(r)$, where $c_{2},c_{3}$ are absolute positive constants. The 
constant
$c_{1}\leq (m+1)\log km+(2k-1)a(q)$ for generic $q$, and
$c_{1}\leq (m+1)\log km+ 3a(q)$ if restriction of functionals $f_{s}$ 
to any complex straight line passing through 0 generates a 
one-dimensional vector space over $\Re$.}

\end{document}